\documentclass[10pt,A4]{article}
\usepackage[top=1in, bottom=1in, left=1.2in, right=1.2in]{geometry}
\usepackage{graphicx}
\usepackage{here}
\usepackage{amssymb}
\usepackage{amsmath}
\usepackage{amsthm}
\usepackage{epstopdf}
\usepackage{here}
\usepackage{setspace}
\usepackage{mathtools}
\usepackage{comment}
\DeclarePairedDelimiter{\abs}{\lvert}{\rvert}
\newtheorem{Theorem}{Theorem}[section]
\newtheorem{Lemma}[Theorem]{Lemma}
\newtheorem{Corollary}[Theorem]{Corollary}
\newtheorem{Property}[Theorem]{Property}
\newcommand{\Proof}{\noindent{\bf Proof}\quad}
\theoremstyle{definition}

\newtheorem{Definition}[Theorem]{Definition}
\newtheorem{Example}[Theorem]{Example}

\title{\textbf{Optimality and Constructions of \\ Spanning Bipartite Block Designs}}
\author{${}^a$Shoko Chisaki,  ${}^b$Ryoh Fuji-Hara and ${}^c$Nobuko Miyamoto\\
\\
${}^a$Department of Information Systems, Osaka Institute of Technology \\[0.1cm]
${}^b$Faculty of Engineering, Information and Systems, University of Tsukuba \\[0.1cm]
${}^c$Department of Information Science, Tokyo University of Science\\
}
\date{}

\begin{document}

\maketitle

\begin{abstract}
  We consider a statistical problem to estimate variables (effects) that are associated with the edges of a complete bipartite graph $K_{v_1, v_2}=(V_1, V_2 \, ; E)$.
  Each data is obtained as a sum of selected effects, a subset of $E$.
  In order to estimate efficiently, we propose a design called Spanning Bipartite Block Design (SBBD).
   For SBBDs such that  the effects are estimable, we proved that the estimators have the same variance (variance balanced).
  If each block (a subgraph of $K_{v_1, v_2}$) of SBBD is a semi-regular or a regular bipartite graph, we show that  the design is A-optimum.
  We also show a construction of SBBD using an ($r,\lambda$)-design and an ordered design.
  A BIBD with prime power blocks gives an A-optimum semi-regular or regular SBBD.
  At last, we mention that this SBBD is able to use for deep learning.
\\[0.3cm]
\textbf{Keywords.} spanning bipartite block design, A-optimum, variance balanced, ($r,\lambda$)-design, balanced incomplete block design, ordered  design, deep learning
\\[0.3cm]
\textbf{AMS classification.} 62K05, 62K10, 05B05
\end{abstract}

\section{Introduction}

\begin{figure}[ht]
\centering
\includegraphics[scale=0.37,clip]{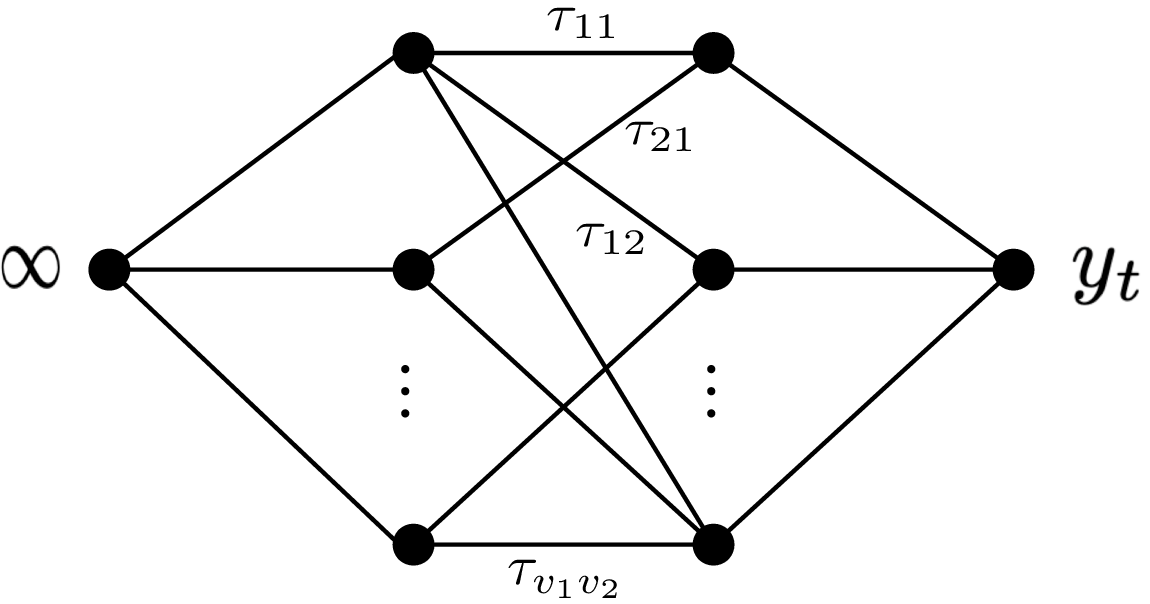}
\caption{Image of bipartite problem}
\label{connect}
\end{figure}

Let $V_1$ and $V_2$ be point sets, and $E$ the set of edges between the $V_1$ and $V_2$, it is a complete bipartite graph, $K_{v_1,v_2} = ( V_1, V_2\, ; E)$, where $\abs{V_1} =v_1, \abs{V_2} =v_2$.
We consider a statistical problem  estimating the variables associated with  $E$ from experimental data.
For example, communication capacities between two sets of cities, traffic volume between two sets of cities, etc (see Fig. \ref{connect}).

Let $\tau_{ij}$  be a variable (or an effect) to estimate corresponding to the edge $(i,j)$, $i\in V_1\, \, j\in V_2$, of the complete bipartite graph $K_{v_1,v_2}$, and $\boldsymbol{\tau}$ be the vector of $\tau_{ij}$ arranged in the following lexicographical order:
\begin{equation}\label{eq:ordert}
\boldsymbol{\tau}=
( \tau_{11}, \tau_{12},\ldots, \tau_{1v_2}\ ;\  \tau_{21},\tau_{22},\ldots, \tau_{2v_2}\ ; \  \cdots \ ;\ \tau_{v_11},\ldots, \tau_{v_1v_2})^t.
\end{equation}
We consider the following statistical model (\ref{SBBDX}) that each data $y_i$ is obtained as a sum of selected effects, i.e. a subset of $\{ \tau_{11}, \tau_{12}, \ldots, \tau_{v_1v_2}  \}$:
\begin{equation}\label{SBBDX}
\begin{split}
&\mathbf{y} =  \mathbf{X} \boldsymbol{\tau} + \boldsymbol{\epsilon} \\
&\sum_{j=1}^{v_2}\tau_{ij} =0 , \ \  1\le i \le v_1\\
&\sum_{i=1}^{v_1}\tau_{i j} =0 ,\ \  1\le j \le v_2,
\end{split}
\end{equation}
where the data vector $\mathbf{y}=(y_1, y_2,\ldots, y_N)^t$ is
assumed that the mean of all data was subtracted
and $\boldsymbol{\epsilon}$ is a vector of random variables of errors with N$(0,\sigma^2)$.
$X$ is a $(0,1)$-matrix of size ($N\times v_1 v_2$) such that $1$ for selected effects and $0$ for other effects in each row.
Our purpose is to estimate all effects with high precision.
The main problem here is how to design the matrix $X$.

\begin{Example} This is an example that each data is obtained as a sum of
selected effects. The effects have  a bipartite graph structure.
Let $v_1=2, v_2= 3$.
\begin{equation*}
    \begin{split}
    y_1 &= \tau_{11} + \tau_{13} +\tau_{22}+ \tau_{23} + \epsilon_1\\
    y_2 &= \tau_{12} + \tau_{13} +\tau_{21}+ \tau_{22} + \epsilon_2\\
     &  \ \ \vdots  \hspace{20mm}  \vdots \\
    y_N &= \tau_{13} + \tau_{21} + \tau_{22}+ \tau_{23} + \epsilon_{N}
    \end{split}
\end{equation*}

\[
\begin{bmatrix}
y_1\\
y_2\\
\vdots\\
y_N
\end{bmatrix}
=
\begin{bmatrix}
1 & 0 & 1 & 0 & 1 & 1\\
0 & 1 & 1 & 1 & 1 & 0\\
 &  & &  \vdots  & & \\
0 & 0 & 1 & 1 & 1 & 1\\
\end{bmatrix}
\cdot
\begin{bmatrix}
\tau_{11}\\
\tau_{12}\\
\tau_{13}\\
\tau_{21}\\
\tau_{22}\\
\tau_{23}
\end{bmatrix}
+
\begin{bmatrix}
\epsilon_1\\
\epsilon_2\\
\vdots\\
\epsilon_N
\end{bmatrix}
\]
\end{Example}

There is a similar model which is a two-way factorial design having a block factor called an incomplete split-block design, see \cite{Ozawa2002}.
The model is the following:
$$\boldsymbol{y}=X_1 \boldsymbol{\alpha}+ X_2 \boldsymbol{\gamma}+X_{12} (\boldsymbol{\alpha \gamma}) + X_3 \boldsymbol{\beta}+\boldsymbol{\epsilon} ,$$
where $X_1,X_2, X_{12}, X_3$ are (0,1)-matrices having exactly one $1$ in each row, $\boldsymbol{\alpha}, \boldsymbol{\gamma}$ are vectors of main effects, $(\boldsymbol{\alpha \gamma})$ is a vector of interaction effects of $\boldsymbol{\alpha}$ and  $\boldsymbol{\gamma}$, $\boldsymbol{\beta}$ is a vector of block effects,  and $\boldsymbol{\epsilon}$ is the error vector.
In our model, there are no main or block effects, and each data is obtained as the sum of interaction effects within a block instead of blocking effects.
Furthermore, we insist on spreading out the interaction effects in each block as much as possible.

In this paper, we propose a new design named  \textit{spanning bipartite block design} for application to the statistical model (\ref{SBBDX}) and
discuss the precision of the  estimators in the designs.
In Section $2$,
a spanning bipartite block design (SBBD) is defined precisely.
 In Section $3$, we discuss optimality of designs.
 Optimality problems of new designs are discussed in \cite{Ozawa2002} and \cite{Lu21}.
We argue the optimality of SBBD and show a design is A-optimum in a certain class.
In Section $4$, a construction of SBBD is shown using an $(r,\lambda)$-design and an ordered design.
Finally, in Section $5$, we mention the relationship between  SBBD and deep learning.

\section{Spanning Bipartite Block Design}
\label{def_X}

Let $\mathcal{B}=\{B_1,B_2,\ldots, B_N \}$ be a collection of subgraphs of $K_{v_1,v_2} = (V_{v_1},V_{v_2}\, ; E)$ called spanning bipartite blocks (SB-blocks).
If $\mathcal{B}$ satisfies the following five conditions, then we call $(K_{v_1,v_2}\, ;\mathcal{B} )$ a \textit{spanning bipartite block design} (SBBD):

\begin{enumerate}%[label=\textbf{\roman*.}]
  \renewcommand{\labelenumi}{(\roman{enumi})}
\item Each subgraph $B_i$ of $\mathcal{B}$ is incident with all points of $V_1$ and  $V_2$. This is called the \textit{spanning condition}.
\item Each edge of $E$ appears in $\mathcal{B}$ exactly $\mu$ times.
\item Any  two edges $e_{ij}, e_{ij'} \in E$ such that $i \in V_1$, $j, j' \in V_2,\, (j \ne j')$ are included together in $\lambda_{12}$ subgraphs in $\mathcal{B}$.
\item Any two edges $e_{ij}, e_{i'j} \in E$ such that  $i, i' \in V_1, \, (i \ne i')$, $j \in V_2$ are included  together  in $\lambda_{21}$ subgraphs in $\mathcal{B}$.
\item Any  two edges $e_{ij}$, $e_{i'j'} \in E$ such that $i, i' \in V_1, \, (i \ne i')$, $j, j'\in V_2,\, (j \ne j')$ are included  together in $\lambda_{22}$ subgraphs in $\mathcal{B}$.
\end{enumerate}

Next, we define a $(0,1)$-matrix $X$, called a \textit{design matrix}, from the SB-blocks.

\begin{itemize}
\item Suppose that the edges $e_{ij}$ of $K_{v_1,v_2}$ are arranged in the  same lexicographical order as Equation (\ref{eq:ordert}).
\begin{equation}\label{eq:ordere}
( e_{11}, e_{12},\ldots, e_{1v_2}\ ;\  e_{21},e_{22},\ldots, e_{2v_2}\ ; \  \cdots \ ;\ e_{v_11},\ldots, e_{v_1v_2}).
\end{equation}
This sequence of edges corresponds to the columns of $X$.
Denote $(e_{ij})$ for the column number which corresponds to the edge $e_{ij}$.
\item Put $X=[x_{k,(e_{ij})}]$, then  $x_{k,(e_{ij})}$  is the element of the $k$-th row and the $(e_{ij})$-th column of $X$.
The design matrix $X$ is defined by the SB-blocks  $B_1,B_2,\ldots,B_N$ as follows:
\begin{gather}
x_{k,(e_{ij})}=
\begin{cases}
1 &\mbox{ if}\ \  e_{ij} \in B_k \\
0 &\mbox{ otherwise}
\end{cases}
\end{gather}
\item $X$ is an $N \times v_1v_2$ matrix.
\end{itemize}
This is a convenient way to represent a $(0,1)$-matrix for checking the conditions.
Let  $X_i$ be an $N \times v_2$ submatrix of $X$  consisting of $v_2$ columns of $X$ corresponding to  $(e_{i1}, e_{i2},\ldots,$ $e_{iv_2})$.
Then the design matrix $X$ is partitioned into $v_1$ submatrices expressed as $X=( X_1 \mid X_2 \mid \cdots$  $\mid X_{v_1})$.
The conditions of a spanning bipartite block design  $(K_{v_1,v_2}\, ;\, \mathcal{B} )$  can be re-expressed using the design matrix $X=( X_1 \mid X_2 \mid \cdots \mid X_{v_1})$ as follows:
\begin{enumerate}
  \renewcommand{\labelenumi}{(\Roman{enumi})}
\item If $\mathcal{B}$ satisfies the condition (i), any row of $X_i$ is not  zero-vector for $1 \leq i \leq v_1$ and
 $\sum_{i=1}^{v_1} X_i$ has no zero element (the spanning condition).
\item If $\mathcal{B}$ satisfies the condition (ii),
all diagonal elements of $X_i^{\,t} X_i$ are $\mu$ for $1 \leq i \leq v_1$.

\item If $\mathcal{B}$ satisfies the condition (iii),
all off-diagonal elements of $X_i^{\,t} X_i$ are $\lambda_{12}$ for $1 \leq i \leq v_1$.

\item If $\mathcal{B}$ satisfies the condition (iv), all diagonal elements of $X_i^{\,t} X_j$ are $\lambda_{21}$ for  $1 \leq i \ne j \leq v_1$.

\item If $\mathcal{B}$ satisfies the condition (v), all off-diagonal elements of $X_i^{\,t} X_j$ are $\lambda_{22}$ for $1 \leq i \ne j \leq v_1$.
\end{enumerate}
$X^{\,t}X$ is called an \textit{information matrix}.
The information matrix of an SBBD is expressed as follows:
\begin{eqnarray}
%\begin{split}
X^t X &=& I_{v_1}\otimes ( X_i^{\, t} X_i) +(J_{v_1}-I_{v_1})\otimes  (X_i^{\,t} X_j )  \nonumber \\
&=&
 I_{v_1} \otimes
\begin{bmatrix}
\mu & \lambda_{12} & \cdots & \lambda_{12}\\
\lambda_{12} & \mu & \cdots & \lambda_{12}\\
\vdots &  \vdots &  \ddots        &  \vdots \\
\lambda_{12} & \lambda_{12} & \cdots & \mu
\end{bmatrix}
+ (J_{v_1}-I_{v_1}) \otimes
\begin{bmatrix}
\lambda_{21} & \lambda_{22} & \cdots & \lambda_{22}\\
\lambda_{22} &\lambda_{21} & \cdots & \lambda_{22}\\
\vdots &  \vdots &  \ddots       &  \vdots \\
\lambda_{22} & \lambda_{22} & \cdots & \lambda_{21}
\end{bmatrix},
%\end{split}
\end{eqnarray}
where  $1 \le i\ne j \le v_1$, and $I_n$ is the identity matrix of size $n$ and $J_n$ is the $(n \times n)$  all-ones matrix.

A matrix expressed by $a I_n + b (J_n-I_n)$ is called \textit{completely symmetric}.
The information matrix above has a double structure of a completely symmetric matrix.  We call the matrix \textit{double completely symmetric}.
A spanning bipartite block design $(K_{v_1,v_2}\, ; \mathcal{B} )$ is denoted as SBBD($v_1,v_2, N\, ; \Lambda$), where $\Lambda= (\mu, \lambda_{12}, \lambda_{21}, \lambda_{22})$.

\begin{Example}\label{Ex:1}

Let
\begin{center}
$ X = [X_1 \mid X_2 \mid X_3 ] = $
 \scalebox{0.75}{$
\left[
\begin{array}{ccc|ccc|ccc}
 0 & 1 & 1 & 1 & 1 & 0 & 1 & 1 & 0 \\
 1 & 0 & 1 & 0 & 1 & 1 & 0 & 1 & 1 \\
 1 & 1 & 0 & 1 & 0 & 1 & 1 & 0 & 1 \\
 0 & 1 & 1 & 0 & 1 & 1 & 1 & 0 & 1 \\
 1 & 0 & 1 & 1 & 0 & 1 & 1 & 1 & 0 \\
 1 & 1 & 0 & 1 & 1 & 0 & 0 & 1 & 1 \\
 0 & 1 & 1 & 1 & 0 & 1 & 0 & 1 & 1 \\
 1 & 0 & 1 & 1 & 1 & 0 & 1 & 0 & 1 \\
 1 & 1 & 0 & 0 & 1 & 1 & 1 & 1 & 0 \\
\end{array}
\right]
$}
\end{center}
be a design matrix of an SBBD.
Then its information matrix is
$$ X^t X=
I_3 \otimes
\left[\begin{array}{ccc}6 & 3 & 3 \\3 & 6 & 3 \\3 & 3 & 6\end{array}\right]
+ (J_3-I_3)\otimes  \left[\begin{array}{ccc}4 & 4 & 4 \\4 & 4 & 4 \\4 & 4 & 4\end{array}\right].
$$
 The design matrix $X$ satisfies the spanning condition since any  row  of  $X_i$ is not the zero-vector and $X_1+X_2+X_3$ does not contain $0$.
So we have an SBBD$(3, 3, 9\, ;\Lambda),$ $\Lambda=(6,3,4,4)$.
\end{Example}

As you can see from the above example, the spanning condition can not be confirmed from the information matrix $X^t X$.
If  $v_1 \ll v_2$, there is a high possibility that the spanning condition is not met.
Such a design in which the spanning condition (I) is not guaranteed is denoted by SBBD$^*$.

%%%%%%%%%%%%%%%%%%%%%%%%%%%%% Optimality
\section{Optimality}

\subsection{Variance balanced}
For a design matrix, we have a statistical problem of whether it is optimum under certain conditions.
There are some criteria for the precision of the estimators (variances of estimators), see \cite{Pukelsheim2006}. Here we use a criterion called A-optimality (or A-criterion).

Let $e_{ij}^{(v)}=(e_1,e_2,\ldots, e_v)$ be a $(0,1)$-vector of length $v$ such that $e_i = 1$, $e_j=-1$ and  $e_k=0, k\ne i,j$.
$(e_{ij}^{(v_1)} \otimes e_{i'j'}^{(v_2)})^t \boldsymbol{\tau}$
for any $1\le i < j \le v_1$ and $1\le i'< j' \le v_2 $ are called \textit{elementary contrasts}.

Suppose an information matrix $X^tX$ has a double structure of completely symmetric matrices $A$ and $B$ of size $(v_1 \times v_1)$ and $(v_2 \times v_2)$, respectively.
Let $\boldsymbol{p}_1, \boldsymbol{p}_2, \ldots, \boldsymbol{p}_{v_1-1}$ be orthonormal  eigenvectors of $A$ orthogonal to $\boldsymbol{1}_{v_1}$, and
also $\boldsymbol{q}_1, \boldsymbol{q}_2, \ldots, \boldsymbol{q}_{v_2-1}$ be similar vectors of  $B$ with size $v_2$, where $\boldsymbol{1}_n$ is the all-one $n$-vector.
A \textit{basic contrast} of $A \otimes B$ is defined by
\begin{equation}\label{eq:cont}
(\boldsymbol{p}_i \otimes \boldsymbol{q}_j )^t \boldsymbol{\tau}.
\end{equation}

Since every $e_{ij}^{(v_1)} \otimes e_{i'j'}^{(v_2)}$ for any $1\le i < j \le v_1$ and $1\le i'< j' \le v_2 $ lies  in the subspace spanned by  $\boldsymbol{p}_1 \otimes \boldsymbol{q}_1 ,\ \boldsymbol{p}_1 \otimes \boldsymbol{q}_2,  \ldots,\  \boldsymbol{p}_{v_1-1}\otimes \boldsymbol{q}_{v_2-1}$,
we here use basic contrasts for the proofs in this section,
although elementary contrasts are commonly used as contrasts.

Let $\theta_{1,1},\theta_{1,2}, \ldots, \theta_{(v_1-1),(v_2-1)}$ be non-zero eigenvalues of $A \otimes B$, and $\boldsymbol{p}_1 \otimes \boldsymbol{q}_1 ,\ \boldsymbol{p}_1 \otimes \boldsymbol{q}_2,  \ldots,\  \boldsymbol{p}_{v_1-1}\otimes \boldsymbol{q}_{v_2-1}$   be orthonormal eigenvectors, that is,
$(\boldsymbol{p}_i \otimes  \boldsymbol{q}_j )^t \boldsymbol{1}_{v_1v_2} = 0$,
$(\boldsymbol{p}_i \otimes  \boldsymbol{q}_j )^t (\boldsymbol{p}_{i} \otimes \boldsymbol{q}_{j}) =1$ and
$(\boldsymbol{p}_i \otimes  \boldsymbol{q}_j )^t (\boldsymbol{p}_{i'} \otimes \boldsymbol{q}_{j'}) =0$,
$i\ne i'$ or $j\ne j'$
then we have
\begin{equation} \label{eq:var}
\mbox{Var}((\boldsymbol{p}_i \otimes \boldsymbol{q}_j )^t \hat{\boldsymbol{\tau}} )
=(\boldsymbol{p}_i \otimes \boldsymbol{q}_j )^t (A\otimes B)^- (\boldsymbol{p}_i \otimes \boldsymbol{q}_j )\sigma^2
= \frac{1}{\theta_{i,j} }\sigma^2,
\end{equation}
where  $(A\otimes B)^-$ is a Moore-Penrose generalized inverse matrix
of $(A\otimes B)$.

\begin{Definition}(A-optimum, \cite{Kiefer74})
Let $\Xi$ be the set of $(N \times v_1 v_2)$ design matrices $X$ with a certain number of  $1'$s.
Assume $X^tX$, $X\in \Xi$, has $(v_1-1)(v_2-1)$ non-zero eigenvalues
$\alpha_1,\alpha_2, \ldots, \alpha_{(v_1-1)(v_2-1)}$.
For a design matrix $X \in \Xi$,
if the sum of $1 / \alpha_i$ is minimum among $\Xi$,
then  $X$ is called \textit{A-optimum relative to $\Xi$}.
\begin{equation}
\mbox{A-optimum:}\  \min_{X \in \Xi} \ \bigg\{ \sum_{1\le i  \le (v_1-1)(v_2-1)} \frac{1}{\alpha_i} \bigg\}
\end{equation}
\end{Definition}

\begin{Definition}(Variance Balanced, \cite{Rao1958})
 If  all variances for the estimators  of basic contrasts are the same, that is, if   $X^t X$ has $(v_1-1)(v_2-1)$ identical  non-zero eigenvalues, then
 the design is called \textit{variance balanced}:
\begin{equation}
\begin{split}
\mathrm{Var} ((\boldsymbol{p}_i \otimes \boldsymbol{q}_j)^t\boldsymbol{\hat{\tau} })
%=(\boldsymbol{p}_i \otimes \boldsymbol{q}_j)^t (X^tX)^- (\boldsymbol{p}_i \otimes \boldsymbol{q}_j) \sigma^2
= \frac{1}{\alpha}\sigma^2
%\mbox{ \ for } i=1,2,\ldots, (v_1-1), j=1,2,\ldots, (v_2-1).
\mbox{\ \ for }  1\le i  \le v_1-1 \mbox{ and } 1\le  j \le v_2-1.
\end{split}
\end{equation}
\end{Definition}
It is known that any  A-optimum design is variance balanced, but the reverse has not been proven.

Consider a statistical model (\ref{SBBDX}) of SBBD.
If we compare to an ordinal two-way factorial model with interactions,
the treatment effects of Equation (\ref{eq:ordert}) have the exact same structure as the  interaction effects whose information matrix is double completely symmetric, see \cite{Ozawa2002}.

\begin{Theorem}\label{Thm:VB}
 SBBD$^*$ is variance balanced whenever all  basic contrasts  of  $\boldsymbol{\tau}$ are estimable.
\end{Theorem}

\begin{Proof}
Let $X^tX$ be a double completely symmetric information matrix from an SBBD$^*$ which can be expressed with four integers, $a, b, c, d$ as follows:
\begin{equation}
X^tX=I_{v_1}\otimes (a I_{v_2} +bJ_{v_2} )+
(J_{v_1} -I_{v_1}) \otimes (c I_{v_2} +d J_{v_2}).
\end{equation}
Consider a Moore-Penrose generalized inverse matrix of $X^tX$
by putting $A_1, A_2, B_1$ and $B_2$ as follows:
\begin{equation*}
A_1=I_{v_1}-\frac{1}{v_1} J_{v_1}, \
A_2= \frac{1}{v_1} J_{v_1}, \
B_1= I_{v_2}-\frac{1}{v_2} J_{v_2}, \
B_2=\frac{1}{v_2} J_{v_2}.
\end{equation*}
 Using the spectral decomposition method,
the information matrix $X^tX$ can be rewritten as follows:
%
% \begin{eqnarray} \label{eq:XtX}
%   \begin{split}
%     X^tX &=
%     (A_1+A_2)\otimes (a (B_1+B_2)+b v_2 B_2) + (v_1 A_2 - (A_1 + A_2))\otimes (c (B_1+B_2)+d v_2 B_2 ) \nonumber \\
%     &= (A_1+A_2)\otimes (a B_1 + (a+bv_2)B_2) +((v_1-1) A_2 - A_1) \otimes (cB_1+(c+dv_2)B_2) \nonumber \\
%     &= (a-c)(A_1\otimes B_1) + (a+bv_2-c-dv_2)(A_1 \otimes B_2) + (a+c(v_1-1))(A_2\otimes B_1)  \nonumber  \\
%     &+ \, (a+bv_2+(v_1-1)(c+dv_2))(A_2 \otimes B_2).
%   \end{split}
% \end{eqnarray}

\footnotesize
\begin{eqnarray} \label{eq:XtX}
    X^tX &=&
    (A_1+A_2)\otimes (a (B_1+B_2)+b v_2 B_2) + (v_1 A_2 - (A_1 + A_2))\otimes (c (B_1+B_2)+d v_2 B_2 ) \nonumber \\
    &=& (A_1+A_2)\otimes (a B_1 + (a+bv_2)B_2) +((v_1-1) A_2 - A_1) \otimes (cB_1+(c+dv_2)B_2) \nonumber \\
    &=& (a-c)(A_1\otimes B_1) + (a+bv_2-c-dv_2)(A_1 \otimes B_2) + (a+c(v_1-1))(A_2\otimes B_1)  \nonumber  \\
    &&+ \, (a+bv_2+(v_1-1)(c+dv_2))(A_2 \otimes B_2).
\end{eqnarray}
\normalsize
Then we can see the eigenvalues $\alpha$, $\beta$, $\gamma$, $\delta$ of $X^tX$:
\begin{equation*}
\alpha =  a-c , \
\beta  = a-c+(b-d)v_2,
\gamma = a+c(v_1-1) , \
\delta = a+bv_2+(v_1-1)(c+dv_2).
\end{equation*}
We are interested in the first term, $\alpha (A_1 \otimes B_1)$, of Equation (\ref{eq:XtX}).
The matrix of the term can be represented as
$$
(A_1 \otimes B_1) = \sum_{\substack{ 1\le i \le v_1-1 \\ 1\le j \le v_2-1}} \theta_{i,j} (\boldsymbol{p}_i\otimes \boldsymbol{q}_j )(\boldsymbol{p}_i\otimes \boldsymbol{q}_j )^t ,
$$
where $\boldsymbol{p}_i$ ($\boldsymbol{q}_j$) are orthonormal eigenvectors of $A_1$ ($B_1$) orthogonal to $\boldsymbol{1}_{v_1}$ ($\boldsymbol{1}_{v_2}$) and $\theta_{i,j}$ is the eigenvalue corresponding to $\boldsymbol{p}_i\otimes \boldsymbol{q}_j$.
From the matrix forms of $A_1$ and $B_1$, non-zero eigenvalues $\theta_{i,j}$ are all $1$.
So, all basic contrasts of  $\boldsymbol{\tau}$,
$(\boldsymbol{p}_i\otimes \boldsymbol{q}_j )^t \boldsymbol{\tau}$ for  $1\le i \le v_1-1$ and  $1\le j \le v_2-1$,  are obtained from the first term.
Therefore, the Moore-Penrose generalized inverse matrix $(X^tX)^-$ including the ordinal inverse matrix  is written as:
\begin{equation}\label{eq:geninv}
(X^tX)^- =\frac{1}{\alpha} (A_1\otimes B_1) +\frac{1}{\beta} (A_1 \otimes B_2)
+ \frac{1}{\gamma} (A_2\otimes B_1) + \frac{1}{\delta}  (A_2 \otimes B_2).
\end{equation}
If $ \alpha, \beta, \gamma, \delta$ are all non-zero,
$(X^tX)^-$ is an inverse matrix,
otherwise, it is a Moore-Penrose generalized inverse matrix.
If at least one of the eigenvalues is 0,  then it is obtained by removing the term from Equation (\ref{eq:geninv}).
Here we put $\alpha \ne 0, \ ( a > c )$ from the assumption that all basic contrasts of effects in the model (\ref{SBBDX}) are estimable,
and  algebraic multiplicity of $\alpha$ is  $(v_1-1)(v_2-1)$.
Using Equation (\ref{eq:var}), it holds that
\begin{equation*}
  \begin{split}
\textrm{Var}((\boldsymbol{p}_{i} \otimes \boldsymbol{q}_{j} )^t \hat{\boldsymbol{\tau}}) &=
(\boldsymbol{p}_i \otimes \boldsymbol{q}_j )^t (X^t X)^- (\boldsymbol{p}_i \otimes \boldsymbol{q}_j )\sigma^2 \\
&=\frac{1}{\alpha}
(\boldsymbol{p}_i \otimes \boldsymbol{q}_j )^t (A_1 \otimes B_1) (\boldsymbol{p}_i \otimes \boldsymbol{q}_j ) \sigma^2 \\
&=\frac{1}{\alpha} \sigma^2
\end{split}
\end{equation*}
for $1\le i  \le v_1-1$ and $1\le  j \le v_2-1$.
Therefore, the theorem is complete.
\qed
\end{Proof}
\vspace{2mm}

The coefficients $a, b, c, d$ in the proof correspond to the parameters of SBBD as follows:
\begin{eqnarray*}
    a = \mu -\lambda_{12} , \ \
    b = \lambda_{12}   ,   \ \
    c = \lambda_{21} - \lambda_{22} , \ \
    d = \lambda_{22} .
\end{eqnarray*}

\begin{Example}
Eigenvalues of $X^tX$ in Example \ref{Ex:1} are
$$ 36,\ 3,\ 3,\ 3,\ 3,\ 3,\ 3,\ 0,\ 0.$$
They include $(3-1)(3-1)=4$ same eigenvalues for  $\alpha$ $( =3 )$.
The parameters of Example \ref{Ex:1} are $\mu=6$,  $\lambda_{12}=3$,  $\lambda_{21}=4$, $\lambda_{22}=4$,  $v_1=3$, $v_2=3$.
Then we have the following four kinds of eigenvalues, $\alpha, \beta, \gamma, \delta$.
\[
\begin{array}{ccll}
\alpha &=& a-c = \mu -\lambda_{12} - \lambda_{21} + \lambda_{22} = 3 &(m_{\alpha}= (v_1-1)(v_2-1)=4 )\\
\beta  &=& a+bv_2-c-dv_2 = 0   &(m_{\beta}=v_1-1=2 ) \\
\gamma &=& a+c(v_1-1) = 3 &(m_{\gamma} = v_2-1 =2)\\
\delta &=& a+bv_2+(v_1-1)(c+dv_2) =36  &(m_{\delta}= 1),
\end{array}
\]
where $m_{\alpha}, m_{\beta}, m_{\gamma}, m_{\delta}$ are multiplicities.
These are consistent with the proof of Theorem \ref{Thm:VB}.
\end{Example}

From the proof, we have the following corollary:
\begin{Corollary} \label{Col:3.3}
If $X$, size ($N \times v_1v_2$), is a design matrix of an SBBD$^*$, where  all contrasts of  effects are estimable, then $X^tX$ has $(v_1-1)(v_2-1)$ non-zero identical eigenvalues.
\end{Corollary}

\subsection{Optimality on semi-regular SBBD}

Let $B=(V_1,V_2\, ; E'), E' \subseteq E $ be a subgraph of the complete bipartite graph $K_{v_1,v_2}=(V_1,V_2\, ; E)$.
If  all points in $V_1$ and  $V_2$ of $B$ have degrees $k_1$ and  $k_2$, respectively, then the subgraph $B$ is called a \textit{semi-regular bipartite subgraph}, and called \textit{regular bipartite subgraph} if $k_1=k_2$.
Let $\mathcal{B}$ be a set of $N$ semi-regular bipartite subgraphs of $K_{v_1,v_2}$,  degrees are  $k_1$ and $k_2$, that is, $v_1k_1=v_2k_2$. Assume $N \ge (v_1-1)(v_2-1)$.
Let $X$ be a design matrix of size $(N \times v_1v_2)$ defined by Section \ref{def_X} with respect to  $\mathcal{B}$.
Now we define a class  $\Omega$ for matrices $X$  that  satisfies  the following conditions:
\begin{itemize}
\item[(C1)] the number of $1$'s in each column of $X$ is $\mu$,
\item[(C2)] $X$ is a design matrix whose rows correspond to blocks of $\mathcal{B}$,
\item[(C3)] the all elementary contrasts of $X^tX$, $(e_{ij}^{(v_1)} \otimes e_{i'j'}^{(v_2)})^t \boldsymbol{\tau}$, are estimable  for $1\le i < j \le v_1$ and $1\le i'< j' \le v_2 $.
\end{itemize}
If $X=( X_1 \mid X_2 \mid \cdots \mid X_{v_1})$ is a matrix in  $\Omega$, then $X$ has the following properties:
\begin{itemize}
    \item any row of $X_i$ has exactly $k_1$ $1$'s for $1 \le i \le v_1$,
    \item any row of $\sum_{i=1}^{v_1} X_i$ is $(k_2,k_2,\ldots,k_2)$.
\end{itemize}
If a matrix $X \in \Omega$ satisfies the conditions 1 to 5 of SBBD,  it is called a \textit{semi-regular SBBD}.
And if all blocks of a semi-regular SBBD are regular bipartite subgraphs then it is called a \textit{regular SBBD}.

Let $\mathbf{z}_1, \mathbf{z}_2, \ldots , \mathbf{z}_{v_1-1}$ be  the orthonormal vectors orthogonal to $\boldsymbol{1}_{v_1}$, and also
$\mathbf{w}_1, \mathbf{w}_2$, $\ldots, \mathbf{w}_{v_2-1}$ be  the orthonormal vectors  orthogonal to $\boldsymbol{1}_{v_2}$ (they are not necessary to be eigenvectors).

\begin{Lemma}
For any $X \in \Omega$,
$X^tX$ has the following eigenvectors whose eigenvalues are all $0$,
\begin{equation}\label{eq:ev1}
\mathbf{z}_1\otimes \boldsymbol{1}_{v_2}, \mathbf{z}_2\otimes \boldsymbol{1}_{v_2}, \ldots, \mathbf{z}_{v_1-1}\otimes \boldsymbol{1}_{v_2},
\end{equation}
and
\begin{equation}\label{eq:ev2}
\boldsymbol{1}_{v_1}\otimes \mathbf{w}_1,  \boldsymbol{1}_{v_1}\otimes \mathbf{w}_2, \ldots , \boldsymbol{1}_{v_1}\otimes \mathbf{w}_{v_2-1},
\end{equation}
\end{Lemma}
\Proof
Let $k = k_1v_1=v_2k_2$.
Since $X^tX\boldsymbol{1}_{v_1v_2} = \mu k \boldsymbol{1}_{v_1v_2}$, $\boldsymbol{1}_{v_1v_2}$ is an eigenvector of $X^tX$ whose eigenvalue is $\mu k$.
Let $X=(X_1 \mid X_2 \mid \cdots \mid X_{v_1})$ and $X_i = [x_{jh}^{(i)} ]$.
The inner product of the $j$-th row of $X$ and $\mathbf{z}_i \otimes\boldsymbol{1}_{v_2}$ is
$$ \sum_{g=1}^{v_1}\sum_{h=1}^{v_2} x_{j h}^{(g)} z_g^{(i)} = k_1\sum_{g=1}^{v_1} z_g^{(i)} = 0
$$
from the semi-regular condition, where $z_g^{(i)}$ is the $g$-th element of $\mathbf{z}_i$.
Similarly, we have
$$ \sum_{g=1}^{v_1}\sum_{h=1}^{v_2} x_{j h}^{(g)} w_h^{(i)} = k_2 \sum_{h=1}^{v_2} w_h^{(i)} = 0 ,
$$
where $w_h^{(i)}$ is the $h$-th element of $\mathbf{w}_i$.
Therefore, the vectors of (\ref{eq:ev1}) and (\ref{eq:ev2}) are eigenvectors of $X^tX$ corresponding to eigenvalue zero.
\qed

Let  $\boldsymbol{u}_1, \boldsymbol{u}_2, \ldots, \boldsymbol{u}_{(v_1-1)(v_2-1)}$ be orthonormal eigenvectors of $X^tX$ which are orthogonal to $\boldsymbol{1}_{v_1v_2}$, (\ref{eq:ev1}) and (\ref{eq:ev2}).
\begin{Lemma}\label{Lem:2}
Every $X$ in $\Omega$ is A-optimum if the eigenvalues corresponding to  $\boldsymbol{u}_i$'s are all equal.
\end{Lemma}

\begin{Proof}
Let $\alpha_i$ be the eigenvalue corresponding to $\boldsymbol{u}_i$.
We have a spectrum decomposition of $X^tX$,

\begin{equation}\label{eq:xx1}
X^tX = \sum_{i=1}^{(v_1-1)(v_2-1)} \alpha_i \boldsymbol{u}_i \boldsymbol{u}_i^t + \mu k \frac{1}{v_1v_2}\boldsymbol{1}_{v_1v_2}(\boldsymbol{1}_{v_1v_2})^t.
\end{equation}
Every $e_{ij}^{(v_1)} \otimes e_{i'j'}^{(v_2)}$ for any $1\le i < j \le v_1$ and $1\le i'< j' \le v_2 $ lies  in the subspace spanned by $\boldsymbol{u}_1, \boldsymbol{u}_2, \ldots, \boldsymbol{u}_{(v_1-1)(v_2-1)}$.
From the condition (C3), $\alpha_i > 0$ for $1 \leq i \leq (v_1-1)(v_2-1)$.
%$i =1,2,\ldots , (v_1-1)(v_2-1)$.
%
Let $\hat{\boldsymbol{\tau}}$ be  the least square estimator of $\boldsymbol{\tau}$.
We have
$$
\textrm{Var}(\boldsymbol{u}_i^t \hat{\boldsymbol{\tau}}) =\boldsymbol{u}_i^t (X^tX)^- \boldsymbol{u}_i \sigma^2 = \frac{1}{\alpha_i}\sigma^2
$$
with  Moore-Penrose generalized inverse of $X^tX$:
$$
(X^tX)^- = \sum_{i=1}^{(v_1-1)(v_2-1)} \frac{1}{\alpha_i} \boldsymbol{u}_i \boldsymbol{u}_i^t + \frac{1}{\mu k } \frac{1}{v_1v_2} \boldsymbol{1}_{v_1v_2}(\boldsymbol{1}_{v_1v_2})^t .
$$
Consider A-criterion
\begin{equation}\label{eq:A-ct}
\sum_{i=1}^{(v_1-1)(v_2-1)} \frac{1}{\alpha_i} \ .
\end{equation}
From Equation (\ref{eq:xx1}), we have
$$
\textrm{tr}(X^tX) = \sum_{i=1}^{(v_1-1)(v_2-1)} \alpha_i  +\mu k = \mu v_1v_2,
$$
that is,
$$ \sum_{i=1}^{(v_1-1)(v_2-1)} \alpha_i =  \mu (v_1v_2 -k ). $$
Since
$$ \frac{(v_1-1)(v_2-1)}{\sum_{i=1}^{(v_1-1)(v_2-1)} \frac{1}{\alpha_i}} \le \frac{1}{(v_1-1)(v_2-1)}  \sum_{i=1}^{(v_1-1)(v_2-1)} \alpha_i ,$$
if  $\alpha_1=\alpha_2=\cdots=\alpha_{(v_1-1)(v_2-1)} =\frac{ \mu (v_1v_2 -k )}{(v_1-1)(v_2-1)}$, then A-criterion (\ref{eq:A-ct}) is minimum.
\qed
\end{Proof}
\vspace{3mm}

From Theorem \ref{Thm:VB} and Corollary \ref{Col:3.3}, an information matrix $X^tX$ of an SBBD has non-zero $(v_1-1)(v_2-1)$ eigenvalues which are equal.
Therefore, we have the following theorem:

\begin{Theorem}
A semi-regular SBBD  is  A-optimum relative to $\Omega$.
\end{Theorem}

\section{Constructions of SBBD}

\subsection{A construction using an $(r,\lambda)$-design and an ordered design}

\begin{Definition}($(r,\lambda)$-design, \cite{Hand2nd})
  Let $V$ be $v$-point set and $\mathcal{B}=\{B_1, B_2,\ldots,$ $B_b\}$ a collection of subsets  (blocks) of $V$.
If $(V, \mathcal{B})$ holds the following conditions, it is called an $(r,\lambda)$-\textit{design}:
  \begin{itemize}
  \item each point of $V$ is contained in exactly $r$ blocks of $\mathcal{B}$,
  \item any two distinct points of $V$ are contained in exactly $\lambda$ blocks of $\mathcal{B}$.
  \end{itemize}
\end{Definition}

  Let $v$ be the number of points and  $b$  the number of blocks, and put $k_i= \abs{B_i}$ as the block size. If the block sizes are constant $k$, then the $(r,\lambda)$-design is called a \textit{balanced incomplete block design} (BIBD) and denoted by a $(v, k, \lambda)$-BIBD.
An $(r,\lambda)$-design is also called a \textit{regular pairwise balanced design}.
It is not hard  to construct because there is no  block size restriction.
Pairwise balanced designs (PBD,  the first condition of $(r,\lambda)$-design  is not required) have been  well studied, and many recursive constructions are known, see \cite{Hand2nd} and \cite{Wilson1974}.
It is not difficult to modify a PBD to be regular.

Let $(V,\mathcal{B})$ be  an $(r,\lambda)$-design, and $H=[x_{ij}]$ be the $(b\times v)$ incidence matrix between  $\mathcal{B}=\{B_1, B_2,\ldots,B_b\}$ and $V=\{a_1,a_2,\ldots, a_v\}$.
\[
x_{ij} =
\begin{cases}
1 & \mbox{ if }  a_j \in B_i \\
0 & \mbox{ otherwise}.
\end{cases}
\]
Then $H^tH$ is expressed as $r I_v + \lambda(J_v -I_v )$.

\begin{Definition}(Ordered design, \cite{Rao1961})
Let $M$ be an $(\eta({n}^2-n )\times s)$-array with entries from  $\mathbf{N}_{n}=\{1,2,\ldots,$ $n\}$.
If $M=[d_{pq}]$ holds the following conditions, it is called an
 \textit{ordered design}, denoted by $OD_{\eta}(s,n)$.

  \begin{enumerate}
\item[(1)] each row of M consists of $s$ distinct elements in $\mathbf{N}_{n}$, where $s \le  n$,
\item[(2)] in any distinct two columns of $M$,  every ordered pair $(x,y)$ of distinct elements in $\mathbf{N}_{n}$ appears on the same rows  exactly $\eta$ times.
  \end{enumerate}
\end{Definition}

In the condition (2), if the  pair $(x,y)$ is not necessary to be distinct, the $(\eta n^2 \times (s+1))$ array is called orthogonal array.
It is well known that there exists an $(\eta (n^2-n) \times s)$ ordered design
if there is an orthogonal $(\eta n^2 \times (s+1))$ array.
We know that there is an orthogonal $(q^2 \times (q+1))$ array, $q$ a prime power, see \cite{Hedyat1999}. That is:
%We have the following well-known result, see :

\begin{Property}\label{ppOD}
For any prime power $q$,  there exists an $OD_1(q, q)$.
\end{Property}

Suppose  $H$ is a $(b \times v)$ incidence matrix of an $(r, \lambda)$-design with $v$ points and $b$ blocks, and let $\mathbf{h}_i$ be the $i$-th  row vector of $H$.
Then, we can obtain a design matrix by arranging the vectors $\mathbf{h}_i$  in according  to  an ordered design $OD_{\eta}(s,b)$  $M= [ d_{pq} ]$ as follows:
\begin{equation*}
  X= [ \mathbf{h}_{d_{pq}} ]=(X_1 \mid X_2 \mid \cdots \mid X_s ).
\end{equation*}
Note that $X$ is of size $(\eta( b^2 -b) \times vs)$, where each
$X_j$  is an $(\eta(b^2-b)\times v)$-submatrix of $X$ in which the row vector of ${H}$  are put  in according to  the $j$-th column of $M$.

\begin{Example}\label{ex:ordered design}
The  ordered design $OD_{1}(3,3)$ represented by the symbols $\{1,2,3\}$ is on the left side of the following matrices.
The design matrix $X$ with the vectors $\boldsymbol{h}_i, i=1,2,3$,  is on the right side.
\[
OD_1(3,3)=
\begin{bmatrix}
1 & 2 & 3 \\
2 & 3 & 1 \\
3 & 1 & 2 \\
1 & 3 & 2 \\
2 & 1 & 3 \\
3 & 2 & 1
\end{bmatrix}
, \hspace{10mm}
X =
\begin{bmatrix}
\mathbf{h}_1 & \mathbf{h}_2 & \mathbf{h}_3 \\
\mathbf{h}_2 & \mathbf{h}_3 & \mathbf{h}_1 \\
\mathbf{h}_3 & \mathbf{h}_1 & \mathbf{h}_2 \\
\mathbf{h}_1 & \mathbf{h}_3 & \mathbf{h}_2 \\
\mathbf{h}_2 & \mathbf{h}_1 & \mathbf{h}_3 \\
\mathbf{h}_3 & \mathbf{h}_2 & \mathbf{h}_1
\end{bmatrix}.
\]
\end{Example}

Regarding the $j$-th row of $X$ as an SB-block $B_j$, $1 \le j \le N$, which is a spanning subgraph of $K_{s,v}$,
we have an SBBD $( K_{s,v}\, ; \mathcal{B})$, where $N=\eta ( b^2 - b)$.

\begin{Lemma}\label{Lm:rl}
Let $H$ be the $(b \times v)$ incidence matrix of an $(r,\lambda)$-design, and
$\mathbf{h}_1, \mathbf{h}_2,\ldots, \mathbf{h}_b$ be the row vectors of $H$.
Then the following equations hold:
\begin{gather}\sum_{i=1}^{b} \mathbf{h}_i ^t  \, \mathbf{h}_i = r  I_v + \lambda (J_v-I_v),
\label{eq:rl1}
\end{gather}
\begin{gather}
\sum_{i=1}^b  \sum_{j=1}^b
 \mathbf{h}_i^t \, \mathbf{h}_j =r^2 J_v.
\label{eq:rl2}
\end{gather}
\end{Lemma}

\begin{Proof}
From the definition of $(r,\lambda)$-design, $ H^t H=r I_v + \lambda (J_v-I_v)$.
Since $\sum_{i=1}^{b} \mathbf{h}_i ^t \mathbf{h}_i = H^t H$, it holds Equation (\ref{eq:rl1}).
 Next,
\begin{eqnarray*}
  \begin{split}
    \sum_{i=1}^b  \sum_{j=1}^b \mathbf{h}_i^t \, \mathbf{h}_j
    =\sum_{i=1}^b  (\sum_{j=1}^b \mathbf{h}_i^t \, \mathbf{h}_j )
    &=\sum_{i=1}^b   \mathbf{h}_i^t  (r,r,,\ldots,r) \\
    &=(r,r,\ldots,r)^t (r,r,,\ldots,r)  \\
    &= r^2 J_v.
  \end{split}
\end{eqnarray*}
\qed
\end{Proof}

\begin{Theorem} \label{thm:2.2}
If there exists an $(r,\lambda)$-design with $b$ blocks and $v$ points, and an ordered design $OD_{\eta} (s, b)$, then
there is a spanning bipartite block design SBBD$^*$($s, v, N\, ; \Lambda$), $N=\eta  ({b}^2-b )$
and $\Lambda =(\mu, \lambda_{12}, \lambda_{21},$ $\lambda_{22})=
( \eta  r (b-1), \eta  \lambda (b-1), \eta   r (r-1),  \eta (r^2-\lambda)\ )$.
\end{Theorem}

\begin{Proof}
Let $H$ be the $(b\times v)$ incidence matrix of an $(r,\lambda)$-design with $b$ blocks and $v$ elements, and
${\mathbf h}_i$ ($1 \leq i \leq b$) be the $i$-th row  of $H$.
$M=[d_{pq}]$ is an ordered design $OD_{\eta} (s,b)$.
Let $X$ be the design matrix by arranging the row vector of $H$ in according to the ordered design $M$.
\[
X= [ {\mathbf h}_{d_{pq}} ]=(X_1 \mid X_2 \mid \cdots \mid X_s ).
\]
First, we compute a  diagonal submatrix $X_q^{\,t}X_q$ of the information matrix $X^t X$.
In $X_q$, $1\le q\le s$,
each vector ${\mathbf h}_i$ appears $\eta (b-1)$ times.
Therefore, from Lemma \ref{Lm:rl}, we have
\[
X_q^{\,t} X_q  = \eta (b-1) \sum_{j=1}^b {\mathbf h}_j^t {\mathbf h}_j
=  \eta (b-1) \cdot ( r I_v + \lambda (J_v-I_v) ).
\]
Next, we have the following off-diagonal submatrices of $X^tX$ for $1 \le q\ne q' \le s$:
\begin{align*}
X_q^{\,t} X_{q'} &= \eta\, (\  \sum_{i=1}^b \sum_{j=1}^b  {\mathbf h}_i^t {\mathbf h}_j - \sum_{i=1}^b {\mathbf h}_i^t{\mathbf h}_i \ )\\
&= \eta ( r^2  J_v -  ( rI_v + \lambda  (J_v-I_v) ) \\
&=  \eta  ( (r^2-r) I_v + (r^2-\lambda) (J_v-I_v) ).
\end{align*}
\qed
\end{Proof}

Let $P$ be a permutation matrix of size $(v\times v)$ (a $(0,1)$-matrix with exactly one $1$ in every row and column).
\begin{Corollary}
If $X=(X_1 \mid X_2 \mid \cdots \mid X_s )$ is an SBBD, then
\[
X^{(P)} =( X_1 P \mid X_2 P \mid \cdots \mid X_s P )
\]
is also an SBBD with the same parameters as $X$.
\end{Corollary}
\begin{Proof}
For $1 \le i, j \le s$,
\[
 (X_i P)^t (X_j P) = P^t (X_i^{\,t} X_j) P = X_i^{\,t} X_j
\]
because every $X_i^t X_j$ is a completely symmetric matrix.
Therefore, it holds that $(X^{(P)} )^t X^{(P)}=X^t X.$
\qed
\end{Proof}

$X^{(P)}$ may include many different rows from $X$ and
 may have some additional linearly independent rows in $X$.
If we have the following combined design matrix
using  permutation matrices $P_1, P_2, ... , P_{u-1}$:

$$
X^{(I,P_1, P_2, ... , P_{u-1})}=
\begin{bmatrix}
X_1 & X_2 & \cdots & X_{s}\\
X_1 P_1 & X_2 P_1 & \cdots & X_{s} P_1\\
\vdots  &  \vdots  &  & \vdots\\
X_1 P_{u-1} & X_2 P_{u-1} & \cdots & X_{s} P_{u-1}
\end{bmatrix}
$$
then it is an SBBD with the parameters
$(s, v, N\, ;\Lambda$), $N=\eta u ({b}^2-b )$
and $\Lambda =
( \eta u r (b-1),  \eta  u\lambda (b-1), \eta  u r (r-1),  \eta u(r^2-\lambda)\, )$.

\begin{Corollary}\label{cor:2.2}
In Theorem \ref{thm:2.2}, if $s > b-r$ then there exists an  SBBD$(s,v,N\, ; \Lambda)$ with the same parameters.
\end{Corollary}
\begin{Proof}
Since any row  $\boldsymbol{h}_i$ is not the  zero vector,  $X_i$ contains no zero row vector.
If $s=b$, each row of $\sum_{q=1}^s X_q$ is exactly equal to $\sum_{i=1}^b \mathbf{h}_i =(r,r,\ldots, r)$. Therefore if  $s > b-r$, any component of $\sum_{q=1}^s X_q$ is greater than or equal to $1$.
\qed
\end{Proof}

\begin{Example}
Consider an $(r,\lambda)$-design
($V = \{0, 1, 2 \}, \mathcal{B} =\{0,1\},\{1,2\}$,\\ $\{0,2\},\{0,1,2\}\}$) with parameters $r=3, \lambda=2,  v=3,  b=4$.
Its incidence matrix  is
\[
H =
\left(
\begin{array}{ccc}
 1 & 1 & 0 \\
 0 & 1 & 1 \\
 1 & 0 & 1 \\
 1 & 1 & 1
\end{array}
\right).
\]
Then the row vectors of $H$ are
\[
{\mathbf h}_1 =
\left(
\begin{array}{ccc}
 1 & 1 & 0
\end{array}
\right), \ \
{\mathbf h}_2 =
\left(
\begin{array}{ccc}
 0 & 1 & 1
\end{array}
\right),\ \
{\mathbf h}_3 =
\left(
\begin{array}{ccc}
 1 & 0 & 1
\end{array}
\right),\ \
{\mathbf h}_4 =
\left(
\begin{array}{ccc}
 1 & 1 & 1
\end{array}
\right).
\]
Now, we have a design matrix $X$ using an ordered design $OD_1(4,4)$,
\begin{center}
$X =$
\scalebox{0.8}{$
\begin{bmatrix}
{\mathbf h}_1 & {\mathbf h}_2 & {\mathbf h}_3 & {\mathbf h}_4\\
{\mathbf h}_2 & {\mathbf h}_1 & {\mathbf h}_4 & {\mathbf h}_3\\
{\mathbf h}_3 & {\mathbf h}_4 & {\mathbf h}_1 & {\mathbf h}_2\\
{\mathbf h}_4 & {\mathbf h}_3 & {\mathbf h}_2 & {\mathbf h}_1\\
{\mathbf h}_1 & {\mathbf h}_3 & {\mathbf h}_4 & {\mathbf h}_2\\
{\mathbf h}_2 & {\mathbf h}_4 & {\mathbf h}_3 & {\mathbf h}_1\\
{\mathbf h}_3 & {\mathbf h}_1 & {\mathbf h}_2 & {\mathbf h}_4\\
{\mathbf h}_4 & {\mathbf h}_2 & {\mathbf h}_1 & {\mathbf h}_3\\
{\mathbf h}_1 & {\mathbf h}_4 & {\mathbf h}_2 & {\mathbf h}_3\\
{\mathbf h}_2 & {\mathbf h}_3 & {\mathbf h}_1 & {\mathbf h}_4\\
{\mathbf h}_3 & {\mathbf h}_2 & {\mathbf h}_4 & {\mathbf h}_1\\
{\mathbf h}_4 & {\mathbf h}_1 & {\mathbf h}_3 & {\mathbf h}_2\\
\end{bmatrix}. $}
\end{center}

It is an SBBD$(4, 3, 12\,; \Lambda)$, $\Lambda = ( \mu, \lambda_{12}, \lambda_{21}, \lambda_{22} ) = ( 9,  6, 6, 7 ),$ with the information matrix
\[
X^t X = I_3 \otimes \left[\begin{array}{ccc}9 & 6 & 6 \\6 & 9 & 6 \\6 & 6 & 9\end{array}\right]
+ (J_3 - I_3) \otimes \left[\begin{array}{ccc}6 & 7 &7 \\7 & 6 & 7  \\7 & 7 & 6\end{array}\right].
\]
\end{Example}

\subsection{Regular and semi-regular SBBD}

In this part, we consider a case that  $(V,\mathcal{B})$ is a $(v,k,\lambda)$-BIBD and introduce the constructions of semi-regular and regular SBBDs.

\begin{Theorem}
If there is a $(v,k,\lambda)$-BIBD with b blocks and an ordered design $OD_\eta(b,b)$ then there exist a semi-regular SBBD$(b,v,b^2-b \, ;\,  \Lambda)$,
where $\Lambda=(\eta r(b-1) , \eta  \lambda(b-1), \eta  r (r-1),  \eta (r^2-\lambda) )$.
\end{Theorem}
\begin{Proof}
Let $H$ be the ($b\times v$) incidence matrix of a  $(v,k,\lambda)$-BIBD with $b$ blocks.
$H$ has $k$ ones in each row and $r=kb/v$ ones in each column.
Consider $X$ is the design matrix constructed from $H$ and $OD_{\eta}(b,b)$.
Each row of $X$ is a SB-block of the SBBD $(K_{b,v}\, ;\, \mathcal{B})$. From the construction above, each SB-block $B\in \mathcal{B}$ consists of  permuted  rows of $H$. If we reconstruct $(b\times v)$-array from an SB-block, the number of $1$s in each row and in each column is exactly the same as $H$.
That is, each SB-block of $X$ is a semi-regular bipartite subgraph of $K_{b,v}$, where the degrees of each point of $V_1$ and $V_2$ are $k$ and $r$, respectively.
\qed
\end{Proof}

A BIBD with $b=v$ is said to be a symmetric BIBD. From Property \ref{ppOD},
an $OD_1(b,b)$ exists if $b$ is a prime power.

\begin{Corollary}
If there is a $(v,k,\lambda)$-BIBD with a prime power number of blocks
then there exists a semi-regular SBBD, and if the $(v,k,\lambda)$-BIBD is symmetric then there exist a regular SBBD.
\end{Corollary}
Table \ref{tbl:primeBIBD} is a list of existing BIBD with a prime power number of blocks less then $100$ selected from the table in \cite{Hand2nd}.
The BIBDs in the list are all  symmetric except one.
We can construct  many A-optimal SBBDs.

\begin{table}[ht]
\begin{center}
    \caption{BIBD with prime power $b$ blocks}

\begin{tabular}{ccccc|c}
$v$ & $b$ & $r$ & $k$ & $\lambda$  & \text{Remark} \\ \hline \hline
 7 & 7 & 3 & 3 & 1 & \text{PG(2,2)} \\
 11 & 11 & 5 & 5 & 2 &  \\
 13 & 13 & 4 & 4 & 1 & \text{PG(2,3)} \\
 19 & 19 & 9 & 9 & 4 &  \\
 23 & 23 & 11 & 11 & 5 & \\
 25 & 25 & 9 & 9 & 3 &  \\
 27 & 27 & 13 & 13 & 6 & \text{27=$3^3$} \\
 31 & 31 & 6 & 6 & 1 & \text{PG(2,5)} \\
 31 & 31 & 10 & 10 & 3 &  \\
 31 & 31 & 15 & 15 & 7 & \text{PG(4,2)} \\
 37 & 37 & 9 & 9 & 2 & \\
 41 & 41 & 16 & 16 & 6 &  \\
 43 & 43 & 21 & 21 & 10 & \\
  47 & 47 & 23 & 23 & 11 & \\
 \hline
\end{tabular}
\quad
\quad
\quad
\begin{tabular}{ccccc|c}
$v$ & $b$ & $r$ & $k$ & $\lambda$  & \text{Remark} \\ \hline \hline
   7 & 49 & 21 & 3 & 7 &  \\
 49 & 49 & 16 & 16 & 5 & \text{49=$7^2$} \\
 59 & 59 & 29 & 29 & 14 &  \\
 61 & 61 & 16 & 16 & 4 & \\
 61 & 61 & 25 & 25 & 10 &  \\
 67 & 67 & 33 & 33 & 16 &  \\
 71 & 71 & 15 & 15 & 3 & \\
 71 & 71 & 21 & 21 & 6 & \\
 71 & 71 & 35 & 35 & 17 & \\
 73 & 73 & 9 & 9 & 1 & \text{PG(2,8)} \\
 79 & 79 & 13 & 13 & 2 &  \\
 79 & 79 & 27 & 27 & 9 & \\
 79 & 79 & 39 & 39 & 19 & \\
  & & & & &\\ \hline
\end{tabular}
\end{center}
  \label{tbl:primeBIBD}
\end{table}

\section{Application to Deep Learning}

\begin{figure}[H]
  \begin{center}
  \includegraphics[clip,scale=0.25]{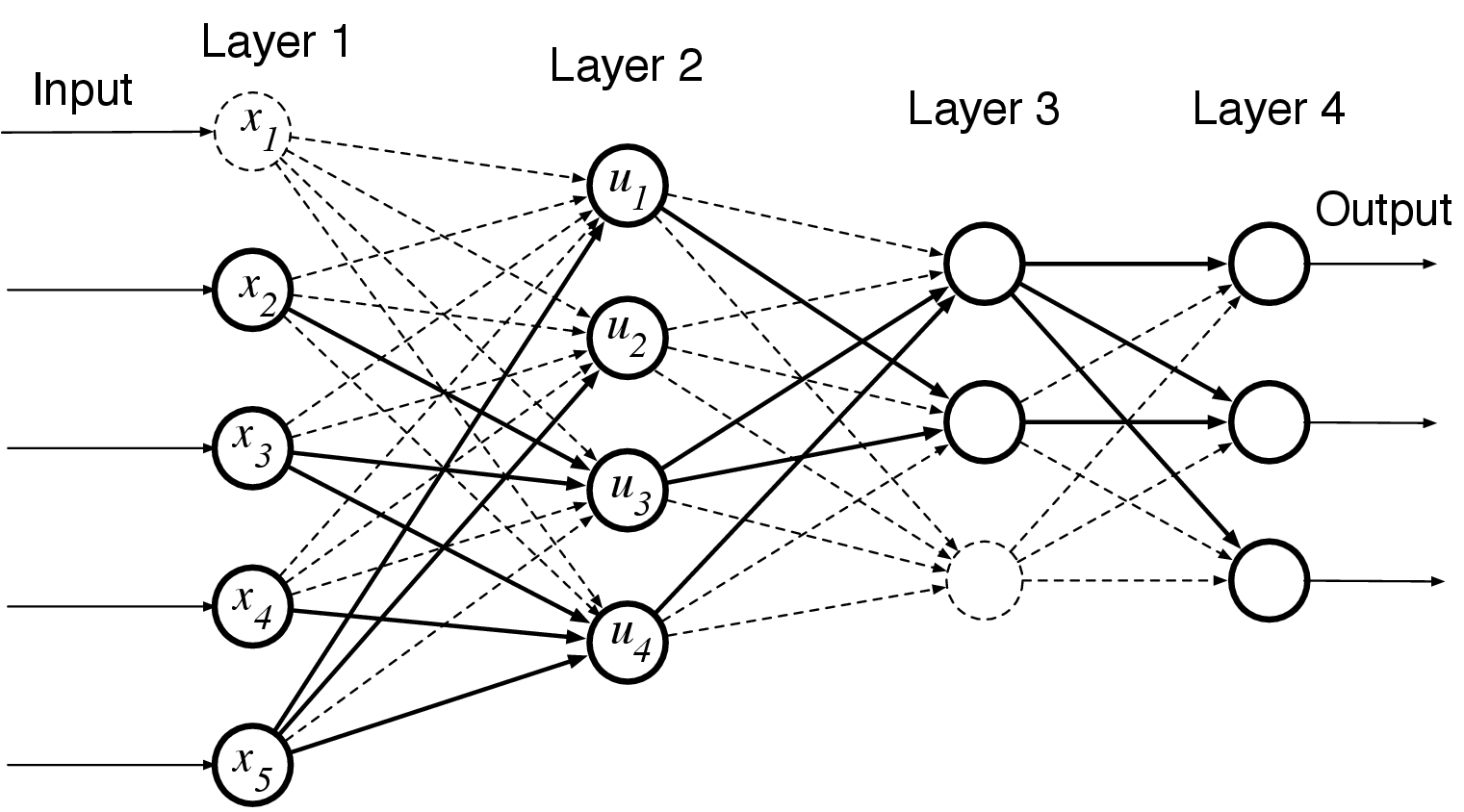}
  \caption{DropConnect}
  \label{fig:DC}
  \end{center}
\end{figure}

Deep learning, in other words, a multi-layer neural network model, is a network consisting of a sequence of point sets (node layers) and complete bipartite graphs (connection layers) between consecutive node layers.
Ignoring the difference between solid and dotted lines, Fig. \ref{fig:DC} gives an example.

\noindent

The weight variables are associated with the edges (connection), and
they are gradually estimated using a large amount of training data that are  pairs of input data $\boldsymbol{x}_i$ and correct answers $\boldsymbol{d}_i, i=1,2,\ldots, N$.
Let  $\boldsymbol{W}$ be the all weight parameters, and  $y(\boldsymbol{x}_i,\boldsymbol{W})$ be output.
$\boldsymbol{W}$ is estimated in the same way as regression so that the following error function is minimized:
    \[
    E(\boldsymbol{W}) = \frac{1}{2}\sum_{i=1}^N \parallel \boldsymbol{d}_i-y(\boldsymbol{x}_i,\boldsymbol{W}) \parallel^2.
    \]
During the learning process, we regularly test using data not used in the learning processes.
At this time, the training data are learning smoothly, but it often gets worse in the test.
This is called \textit{overlearning} or \textit{overfitting}.
It is known that similar overfitting occurs in regression when the model has a large  number of parameters  to be estimated.
In deep learning, \cite{JMLR2014}
proposed a method called \textit{Dropout} in 2014 as a way to deal with overfitting.
This is a method of randomly invalidating the points of each node layer and joining only the valid points by a complete bipartite subgraph.
This is a kind of so-called \textit{sparsification}.
Regarding this method, Chisaki et al. 2020 and 2021 \cite{Chisa1, Chisa2}
have proposed a method applying the combinatorial design theory.

In 2013, \cite{wan2013}
proposed another method called \textit{DropConnect},
it is a method of sparsification by  randomly selecting some  edges in a connection layer instead of  points of a node layer, for an example, solid lines in Fig. \ref{fig:DC}.
We propose to sparsify the edges in connection layers using SBBD, not by random selection.
In the DropConnect method, a node without an incoming connection can occur.
From the spanning condition of SBBD, we can sparsify  independently for each connection layer without occurring of such nodes.
We expect to have a balanced sparsifying system for a multi-layer neural network which have statistically high precision for weight parameter estimations.

\subsection*{Acknowledgments}
I would like to express our sincere gratitude to Professor Shinji Kuriki. He provided very important comments and advice during the writing of this paper. In particular, we received great help from him for the optimality proof.
 This work was supported in part by JSPS KAKENHI Grant Number JP19K11866.

%%%%%%%%%%%%%%%%%%%%%%
% \bibliographystyle{abbrv}
% \bibliography{dropout}
\bibliographystyle{amsplain}
\bibliography{sn-bibliography}

\end{document}